\documentclass[12pt]{article}
\usepackage{amsmath,amsfonts,amssymb,amscd}
\newtheorem{theorem}{Theorem}

\newtheorem{lemma}[theorem]{Lemma}

\newtheorem{definition}{Definition}

\title{Motivic Milnor fibre for nondegenerate function germs on toric singularities}
\author{J.H.M. Steenbrink\\
IMAPP, Radboud University Nijmegen\\ email: J.Steenbrink@math.ru.nl}
\date{25 October 2013}

\begin{document}
\maketitle
\abstract{We study function germs on toric varieties which are nondegenerate for their Newton diagram. We express their motivic Milnor fibre  in terms of their Newton diagram. We extend a formula for the motivic nearby fibre to the case of a toroidal degeneration. We illustrate this by some examples.}

\section{Introduction}
In the calculation of invariants of isolated hypersurface singularities, Newton diagram methods have always been an important tool. These methods are closely connected with the theory of torus embeddings. They were first used for the computation of  resolutions of cusp singularities by Ehlers \cite{Ehlers (1975)}. Kouchnirenko in \cite{Kouchnirenko (1976)} uses them to compute the Milnor number for nondegenerate functions, and Varchenko in \cite{Varchenko (1976)} to compute the zeta function of their monodromy. A conjecture for the Hodge spectrum of nondegenerate hypersurface singularities in terms of their Newton diagram  was formulated in \cite{Steenbrink (1977A)} and  proved by Danilov in \cite{Danilov (1979)}. It was reformulated by Saito  in terms of the Newton filtration on the jacobian module and proved in \cite{Saito (1988)}. 

In this paper we deal with the motivic Milnor fibre of nondegenerate function germs. We consider these as germs defined on a toric variety at a zero-dimensional torus orbit, thereby widening the range of examples of singularities where Newton diagram methods can be applied. We use the Newton diagram to define a toric modification of the ambient space which is a toric variety and such that the zero fiber of our germ defines a toroidal divisor on it. The proof of our formula relies then on the extension of a formula for the motivic nearby fiber from the normal crossing case to the toroidal case. 

The importance of the motivic Milnor fibre lies in the fact that several additive invariants like Euler characteristic, Zeta function of monodromy and certain invariants connected with mixed Hodge structures can be read from this invariant. 

 Each of the terms occurring in our formula for the motivic Milnor fibre is a hypersurface in a torus given by an equation $g(z)=1$ with $g$ weighted homogeneous; we will calculate the Hodge spectrum of such hypersurfaces in a subsequent paper with Sander Rieken. 

More friends I made during my mathematical career who deserve a birthday present than I am able to produce papers of the quality they deserve. Therefore I dedicate this paper to the sixtiest birthdays of Sabir Gusein-Zade, Wolfgang Ebeling and Alexandru Dimca. I hope they do not mind to share this present.  

\section{Motivic nearby fibre}
We briefly recall the notions of motivic nearby fibre and motivic Milnor fibre as defined by Denef and Loeser \cite{Denef and Loeser (2001)}. We prove a  formula for them  in a toroidal setting.

 In this section, $k$ is an algebraically closed field of characteristic zero. 

\subsection{Grothendieck groups of varieties}
Let us recall the Grothendieck group $K_0(\mathrm{Var}_k)$ of varieties over $k$. It is defined as an  abelian group  by the following generators and relations. Generators are the  isomorphism classes $[X]$ of algebraic $k$-varieties and relations are of the form $[X] = [Y]+[X \setminus Y]$ for any pair $Y\subset X$ where $Y$ is a closed subvariety. By defining the product of $[Z]$ and $[W]$ to be $[Z \times W]$ we obtain a commutative ring with unit $[\mathrm{Spec}(k)]$. 

Let $X$ be a smooth projective variety and $Y \subset X$ a smooth closed subvariety. Let $\pi:X'\to X$ be the blowing-up with center $Y$ and let $Y'= \pi^{-1}(Y)$. Then $\pi$ defines an isomorphism between $X' \setminus Y'$ and $X \setminus Y$, hence 
\begin{equation} \label{blowup}
[X'] - [Y'] = [X] - [Y] \mbox{ in } K_0(\mathrm{Var}_k)\ .
\end{equation}
Bittner \cite{Bittner (2004)} showed that this type of  relation suffices to describe the group $ K_0(\mathrm{Var}_k)$:
\begin{theorem}
The group $ K_0(\mathrm{Var}_k)$ is isomorphic to the abelian group with generators $[X]$ where $X$ is smooth projective and relations (\ref{blowup}).
\end{theorem}

If $S$ is a $k$-variety, we have a relative Grothendieck group $K_0(\mathrm{Var}_S)$ with generators $[X]_S$ the $S$-isomorphism classes of varieties over $S$ and again  relations $[X] = [Y]+[X \setminus Y]$ for any pair $Y\subset X$ where $Y$ is a closed subvariety. Fibre product over $S$ equips $K_0(\mathrm{Var}_S)$ with the structure of a ring. We let $\mathbb{L} = [\mathbb{A}^1\times S \to S] \in K_0(\mathrm{Var}_S)$, the map to $S$ being the projection to the second factor, and $\mathcal{M}_S = K_0(\mathrm{Var}_S)[\mathbb{L}^{-1}]$. 

We also need equivariant versions of these constructions. For $n \in \mathbb{N}$ we let $\mu_n$ denote the group of $n$-th roots of unity in $k$. By mapping $\mu_{nd}$ to $\mu_n$ by $x \mapsto x^d$ we obtain a projective system and we let $\hat{\mu} = \varprojlim \mu_n$. 

A \emph{good $\hat{\mu}$-action} on an $S$-variety $X$ is given by an action of $\mu_n$ on $X$ by $S$-morphisms for some $n$ with the property that each orbit is contained in an affine open subset. The group $K_0(\mathrm{Var}^{\hat{\mu}}_S)$ has generators $[X,\hat{\mu}]$ where $X$ is an $S$-variety with good $\hat{\mu}$-action and relations $[X,\hat{\mu}] = [Y,\hat{\mu}]+[X \setminus Y,\hat{\mu}]$ and $[X\times V,\hat{\mu}] = [X,\hat{\mu}]\cdot \mathbb{L}^m$ when $V$ is an $m$-dimensional affine space with any good $\hat{\mu}$-action. Finally $\mathcal{M}_S^{\hat{\mu}} = K_0(\mathrm{Var}^{\hat{\mu}}_S)[\mathbb{L}^{-1}]$.

If $f:S \to T$ is a $k$-morphism, then every $S$-variety $h:X \to S$ becomes also a $T$-variety $fh:X \to T$. This defines group homomorphisms 
$f_!:K_0(\mathrm{Var}_S) \to K_0(\mathrm{Var}_T)$ and $f_!:\mathcal{M}_S \to \mathcal{M}_T$. In the equivariant setting we obtain 
\begin{equation}
f_!:\mathcal{M}_S^{\hat{\mu}} \to \mathcal{M}_T^{\hat{\mu}}.
\end{equation}
On the other hand, if $Y$ is a $T$-variety, then $Y\times_T S$ is an $S$-variety. This defines algebra homomorphisms 
\begin{equation}
f^\ast:  \mathcal{M}_T^{\hat{\mu}} \to \mathcal{M}_S^{\hat{\mu}}.
\end{equation}
\subsection{Nearby and Milnor fibre}
Let $X$ be a smooth connected quasi-projective variety over $k$ and let $f:X \to k$ be a non-constant regular function. Assume that $E = f^ {-1}(0)$ is a divisor with strict normal crossings on $X$. Let $E_i,\ i \in I$, be the irreducible components of $E$ and let $e_i$ denote the multiplicity of $f$ along $E_i$.  

Choose a common multiple $e$ of all $e_i,\ i \in I$. Let $\tau: k \to k$ be defined by $\tau(z) = z^e$, and let $\tau^\ast X$ denote the pull-back of $X$ via $\tau$. Finally let $\tilde{X}$ be the normalization of $\tau^\ast X$. Then we have a diagram
\begin{equation}
\begin{array}{ccccc}
D &\hookrightarrow  &\tilde{X} & \stackrel{\rho}{\rightarrow} & X  \\
\downarrow & & \tilde{f} \downarrow & & f \downarrow \\
0 & \in &k & \stackrel{\tau}{\rightarrow} & k
\end{array}
\end{equation}
Let $\zeta \in \mu_e$. The covering transformation $z \mapsto \zeta z$ of $\tau$ extends to an automorphism $\gamma(\zeta)$ of order $e$ of $\tilde{X}$ and in this way we obtain a good $\hat{\mu}$-action on $\tilde{X}$ and $D$. 

We have  decompositions of the zero fibres $E$ and $D$ into nonsingular locally closed subsets as follows. 
Let $J \subset I$ be a non-empty subset. We define $E_J = \bigcap_{i \in J} E_i$ and $E_J^0 = E_J \setminus \bigcup_{i \not\in J} E_i$. Moreover we let $D_J^0 = \rho^{-1}(E_J^0)$, and $\gamma_J$ the $\hat{\mu}$-action on $D_J^0$ induced by $\gamma$. Then $\rho_J:D^0_J \to E_J^0$ is an etale covering of degree equal to the greatest common divisor of the multiplicities $e_i, i \in J$. 

The \emph{motivic nearby fibre} of $f$ is defined as 
\begin{equation} \label{nearby}
\psi_f = \sum_J [(D_J^0,\gamma_J)](1-\mathbb{L})^{\sharp J -1} \in K_0(\mathrm{Var}_E^{\hat{\mu}} )
\end{equation}
where the sum runs over the non-empty subsets $J$ of $I$. See \cite{Looijenga (2002)}. 

\paragraph{Remarks}
\begin{enumerate}
\item The pairs $(D_J,\gamma_J)$ do not depend on the choice of the integer $e$ (as long as $e_i | e$ for each $i \in I$).
\item The element $\psi_f$ does not change when the zero fiber $E$ is modified by blowing-up. 
\end{enumerate} 

In \cite{Denef and Loeser (2001)} the motivic nearby fibre $\psi_f$  has been defined for any $f:X \to k$ with $X$ nonsingular via the theory of arc spaces. Their definition leads to formula (\ref{nearby}) when the zero fibre $X_0$ of $f$ has strict normal crossings. Under this hypothesis  $\psi_f$  is characterized by the following properties:
\begin{enumerate}
\item if $\pi:X'\to X$ is a proper birational map which induces an isomorphism $X'\setminus X'_0 \to X \setminus X_0$ and $\pi_0:X'_0 \to X_0$ is the induced map, then $\psi_f = \pi_{0!} \psi_{f\pi}$;
\item if $X_0$ is a divisor with strict normal crossings, then $\psi_f$ is given by  formula (\ref{nearby}).
\end{enumerate} 
These two properties enable one to extend the definition of $\psi_f \in K_0(\mathrm{Var}_{X_0}^{\hat{\mu}} )$ to the case where $X$ is singular but $\mathrm{Sing}(X) \subset X_0$: one chooses an embedded resolution $\pi: (X',E)\to (X,X_0)$ such that $E$ is a divisor with strict normal crossings and defines  $\psi_f = \pi_{0!} \psi_{f\pi}$.

In the sequel we will need a definition of $\psi_f$ in the general setting where $X$ may have singularities not contained in $X_0$. Here we use Bittner's \emph{ nearby cycle functor} 
\begin{equation}
\Psi_f: \mathcal{M}_X \to \mathcal{M}_{X_0}^{\hat{\mu}}
\end{equation}
as defined in \cite{Bittner (2005)}. Bittner shows that $\mathcal{M}_X$ is generated by $[Y]_X$ with $Y$ proper over $X$ and nonsingular, and for such $[Y]$ she defines  
\begin{equation}
\Psi_f([Y]) = \pi_{0!}\psi_{f\pi},
\end{equation}
where $\pi:Y \to X$ is the given morphism and $\pi_0: Y_0 \to X_0$ its restriction to the zero fibres.
This formula respects the relations between generators by  \cite[Claim 8.2]{Bittner (2005)}.  Hence we obtain an $\mathcal{M}_k$-linear map $\Psi_f: \mathcal{M}_X \to \mathcal{M}_{X_0}^{\hat{\mu}}$ and may define
\begin{equation}
\psi_f := \Psi_f([X]).
\end{equation}
Following \cite{Denef and Loeser (2001)} for $x \in X_0$ we define the \emph{motivic Milnor fibre of $f$ at $x$} by 
\begin{equation}
\psi_{f,x} = i_x^\ast\psi_f
\end{equation}
where $i_x:\{x\} \hookrightarrow X_0$ is the inclusion. 
\subsection{Toroidal embeddings} \label{torus}
We first recall the basic notions about toric varieties. Details can be found in \cite{Fulton (1993)}, \cite{Danilov (1978)} or \cite{Oda 1988}. 

An $n$-dimensional \emph{toric variety} is a normal variety $X$ which contains an algebraic torus $T \simeq (\mathbb{G}_m)^n$ as a dense Zariski open subset such that the action of $T$ on itself by translation extends to an action of $T$ on $X$ with finitely many orbits. 

We let $M$ denote the character group of $T$; it is free abelian of rank $n$. The dual group $N = \mathrm{Hom}(M,\mathbb{Z})$ can be identified with the group of one-parameter subgroups of $T$. We define $M_\mathbb{R} = M \otimes \mathbb{R}$ and $N_\mathbb{R} = N \otimes \mathbb{R}$. 

Each $T$-invariant affine open subset $U$ of $X$ corresponds to a cone in $N_\mathbb{R}$ which is the convex hull of all one-parameter subgroups of $T$ which extend to a morphism $\mathbb{A}^1_k \to U$. Conversely, the affine variety $U$ is determined by its cone $\sigma \subset N_\mathbb{R}$ as follows. Let $\sigma^\vee = \{m \in M_\mathbb{R} \mid \forall p \in \sigma: \langle m,p \rangle \geq 0 \}$. Then 
\begin{equation}
U = X_\sigma := \mathrm{Spec}\ k[M \cap \sigma^\vee].
\end{equation}
The (finite) collection $\Sigma$ of these cones is called the \emph{fan} in $N_\mathbb{R}$ associated to $X$, and $X$ can be recovered from its fan. The variety $X$ is complete if $\bigcup_{\sigma \in \Sigma} \sigma = N_\mathbb{R}$.  
\begin{definition}
A \emph{toroidal embedding (without self-intersection)}  is an open embedding $U \hookrightarrow X$ of varieties with the following property: for every $x \in X$ there exists an open neighborhood $U_x$ and an etale morphism $\phi:U_x \to X_{\sigma_x}$ to a toric variety $X_{\sigma_x}$ with torus $T_x$ such that $U_x \cap U = \phi^{-1}(T_x)$. Such a morphism is called a \emph{chart at $x$}.
\end{definition}
If $U \hookrightarrow X$ is a toroidal embedding such that $X$ is nonsingular, then $X \setminus U$ is a divisor on $X$ with strict normal crossings. 
For every toroidal embedding $U \hookrightarrow X$ we have a canonical stratification with the following properties: the strata are connected locally closed subvarieties of $X$ and for every chart $\phi:U_x \to X_{\sigma_x}$ the intersection of the strata with $U_x$ are the preimages of the torus orbits in $X_{\sigma_x}$ under $\phi$. 
\begin{definition}
A \emph{proper modification} of toroidal embeddings $[U\hookrightarrow X'] \to [U \hookrightarrow X]$ is a proper morphism $\pi:X'\to X$ which is the identity on $U$ (hence birational) and which maps strata of $X'$ onto strata of $X$. 
\end{definition}
\begin{theorem}
For every toroidal embedding $U \hookrightarrow X$ there exists a proper modification of toroidal embeddings $[U\hookrightarrow X'] \to [U \hookrightarrow X]$ such that $X'$ is nonsingular.
\end{theorem} 
This follows from \cite[Theorem 11*]{Kempf et al (1974)}.

\subsection{Motivic nearby fibre in toroidal setting}
We first explain what is meant by ``toroidal setting''. We consider a toroidal embedding $U \subset X$ and a regular function $f:X \to k$ with the property that $U\cap X_0 = \emptyset$. As an important example we have the case where $X$ is nonsingular and $X_0$ is a divisor with strict normal crossings on $X$; here we take $U = X_0$. In general, $\mathrm{Sing}(X)$ will  not be contained in $X_0$; then $U$ has to be strictly smaller than $X \setminus X_0$. 

We are going to show that in the toroidal setting an explicit formula for $\psi_f$ analogous to (\ref{nearby}) is valid. The strata $E_J^0$ in that formula will be replaced by the strata of the toroidal boundary $X \setminus U$ which are contained in $X_0$. 

We start our computation by analyzing the strata of the canonical toroidal stratification of $X$. Let $S$ denote the set of strata, and $S_h$ the subset of $S$ consisting of all strata which are not contained in $X_0$. For $s \in S_h$ let $\bar{s}$ be its closure in $X$ and $i_s: \bar{s} \hookrightarrow X$ the inclusion. 

We have 
\begin{equation}
[X] = [X_0] + \sum_{s \in S_h} [s] 
\end{equation}
in $\mathcal{M}_X$ and hence 
\begin{equation}
\psi_f := \Psi_f([X]) = \sum_{s \in S_h} \Psi_f([s])
\end{equation}
in $\mathcal{M}_{X_0}^{\hat{\mu}}$, because $\Psi_f$ is $\mathcal{M}_k$-linear and $\Psi_f([X_0])=0$ by \cite[Properties 8.4]{Bittner (2005)}. As $i_s$ is a closed embedding, we also have 
\begin{equation}
 \Psi_f([s]) = \Psi_f i_{s!}[s] = i_{s_0!}\Psi_{fi_s}[s]
\end{equation}
where $i_{s_0}$ is the inclusion of $\bar{s}\cap X_0$ in $X_0$. Note that $s \subset \bar{s}$ is also a toroidal embedding. So we have reduced the computation to the case of $\Psi_f([U])$.

\begin{lemma}
Let $U \subset X$ be a toroidal embedding and $f:X \to k$ a regular function such that $U \cap X_0 = \emptyset$. Let $\pi:X'\to X$ be a toroidal modification such that $X'$ is smooth. Let $f'=f\pi$. Then 
\begin{equation}
\Psi_f([U]) =\pi_{0!} \Psi_{f'}([U]).
\end{equation}
\end{lemma}
{\em Proof }
By \cite[Properties 8.4]{Bittner (2005)} we have $\Psi_f\pi_! = \pi_{0!}\Psi_{f\pi}$ because $\pi$ is proper. Moreover $\pi_!([U]_{X'}) = [U]_X$. 

\vspace{4mm}\noindent
Hence to compute $\Psi_f([U])$ we may assume that $X$ is smooth and $X \setminus U$ is a divisor with strict normal crossings. 
\begin{theorem}\label{psiu}
Let $f:X \to k$ be a regular function on a smooth variety $X$. Let $Y \subset X$ be a closed subset such that $Y \cup X_0$ is a divisor with strict normal crossings on $X$. Let $\tilde{U} = X \setminus Y \stackrel{j}{\hookrightarrow} X$ and $\tilde{f} = fj$. Then 
\begin{equation}
\Psi_f([U]) = j_{0!}\psi_{\tilde{f}}.
\end{equation}
Moreover, $\psi_{\tilde{f}}$ is given by formula (\ref{nearby}) with $E = X_0 \setminus Y$. 
\end{theorem}
{\em Proof }
We have $[U] =  \sum_{J \subset I} (-1)^{|J|} [Y_J]$ in $K_0(\mathrm{Var}_X)$, where $I$ is the set of irreducible components of $Y$ and $Y_J = \bigcap_{j\in J} Y_j$ (so $X = Y_\emptyset$). The inclusion $i_J: Y_J \to X$ is proper so $\Psi_f([Y_j]) = (i_J)_{0!}\psi_{fi_J}$, which is computed using formula (\ref{nearby}) with $E = Y_J \cap X_0$. Putting all these terms together we obtain the result.

\begin{theorem}
Let $\pi:X'\to X$ be a proper equivariant modification of $n$-dimensional toroidal embeddings with good $\hat{\mu}$-action. Let $s$ be a stratum of $X$ and let $I'_s$ denote the set of strata of $\pi^{-1}(s)$. Let $c(t) = n - \dim t$ for any stratum $t$. Then 
\begin{equation}\label{fibre}
(1-\mathbb{L})^{c(s)}[s] = \pi_! \sum_{t \in I'_s} (1-\mathbb{L})^{c(t)}[t] \mbox{ in } K_0(\mathrm{Var}^{\hat{\mu}}_X).
\end{equation}
\end{theorem}
{\em Proof }
Each stratum $t\in I'_s$ is a trivial $\mathbb{G}_m^k$-bundle over $s$ with $k = c(s)-c(t)$. Hence 
\begin{equation}
\pi_![t] = (\mathbb{L} - 1)^{c(s)-c(t)}[s]
\end{equation}
so
\begin{equation}
\pi_!(1-\mathbb{L})^{c(t)}[t] = (-1)^{c(s)-c(t)}(1-\mathbb{L})^{c(s)}.
\end{equation}
We claim that $\sum_{t\in I'_s} (-1)^{c(t)} = (-1)^{c(s)}$. By restricting to a chart intersecting $s$  we obtain the situation where all strata are pulled back from a toric variety. So for the proof of (\ref{fibre}) we may assume that $\pi$ is a modification of torus embeddings. Let $\Sigma$ be the fan corresponding to $X$ and $\Sigma'$ the subdivision of $\Sigma$ corresponding to $\pi$. For $\sigma \in \Sigma$ let $\sigma^0$ denote the relative interior of $\sigma$; it is homeomorphic to a cell of dimension $c(s)$, where $s$ is unique closed torus orbit in $X_\sigma$. The orbits $t\subset  \pi^{-1}(s)$ correspond to the cones $\tau \in \Sigma'$ for which $\tau^0 \subset \sigma^0$ and the interiors of htese cones form a cell decomposition of $\sigma^0$. Then 
\begin{equation}
\sum_{t \in I_s} (-1)^{c(t)} = (-1)^{c(s)}
\end{equation}
is a consequence of the additivity 
 the Euler characteristic with compact supports, because an open cell of dimension k has $\chi_c = (-1)^k$. 

\section{Nondegenerate Laurent polynomials}
In this section we consider function germs on toric singularities. We compute their motivic Milnor fibre. 
\subsection{Toric and toroidal singularities} As before, $k$ is an algebraically closed field of characteristic zero and $T$ is an $n$-dimensional algebraic torus over $k$ with character group $M$ as in Sect.~\ref{torus}. We keep the notations of that section. 
\begin{definition} A \emph{toric singularity} is a germ $(X,x)$ where $X$ is a toric variety and $x$ is a fixed point of the torus action on $X$. A \emph{toroidal singularity} is a germ which is analytically isomorphic to a toric singularity. 
\end{definition}
In studying a toric singularity $(X,x)$ we may suppose that the representative $X$ is affine (else replace $X$ by the star of $x$). Then $X = X_\sigma = \mathrm{Spec}\ k[M \cap \sigma^\vee]$ for a strictly convex rational polyhedral cone $\sigma$  in $N_\mathbb{R}$. We let $A_\sigma =  k[M \cap \sigma^\vee]$.

\subsection{Newton polyhedron and nondegeneracy} 
 A regular $k$-valued function on $T$ is called a \emph{Laurent polynomial}. Any character $m$ of the torus can be considered as a regular function $e(m): T \to k$, and the $e(m)$ with $m \in M$ form a k-basis of the coordinate ring $k[T]$. 
\begin{definition}
Given a Laurent polynomial $f = \sum_{m \in M} a_me(m)$ we define its \emph{support} by 
$\mathrm{supp}(f) := \{m \in M \mid a_m \neq 0\}$. For any $V \subset M$ we let $f^V := \sum_{m \in M\cap V} a_me(m)$
\end{definition}
\begin{definition}
Consider a toric singularity $(X_\sigma,x)$  and the corresponding $k$-algebra $A_\sigma \subset k[T]$. The \emph{Newton polyhedron} $\Delta $ of $f$ is the convex hull of $\bigcup_{m\in \mathrm{supp}(f)} (m+\sigma^\vee)$ in $M_\mathbb{R}$ .

The function $f$ is called \emph{convenient} if $\sigma^\vee \setminus \Delta$ is bounded. Equivalently: $\mathrm{supp}(f)$ has non-empty intersection with each one-dimensional face of $\sigma^\vee$. 

The function $f$ is called \emph{Newton non-degenerate} if for each  compact face $\Gamma$ of $\Delta$ the functions $x_1\frac{\partial f^\Gamma}{\partial x_1},\ldots,x_n\frac{\partial f^\Gamma}{\partial x_n}$ have no common zero in the torus $T$.
\end{definition}
This concept is similar to $\Delta$-regularity as in \cite{Batyrev (1993)}. 

 Let $P$ be an integral convex polytope in $M_\mathbb{R}$, i.e. the convex hull of a finite subset of $M$. We let $L(P)$ denote the $k$-linear span of all monomials $e(m)$ with $m \in M \cap P$. If $P$ has dimension $n$, these monomials embed the torus $T$ in projective space, and the associated toric variety $\mathbb{P}(P)$ is the closure of the image $T(P)$ of $T$ under this mapping. It is also equal to $\mathrm{Proj}\ S(P)$ where $S(P) = \bigoplus_{k \in \mathbb{N}} L(kP)$. 

The space $L(P)$ is in a natural way the space of global sections in a very ample line bundle $\mathcal{O}(1)$ on $\mathbb{P}(P)$. Any nonzero element $g \in L(P)$ therefore determines a hyperplane section $\bar{Z}_g$ of $\mathbb{P}(P)$. Such $g$ is called $P$-\emph{regular} if $P$ is the convex hull of $\mathrm{supp}(g)$ and for each  face $Q$ of $P$ the polynomial equations
$$
g^Q = g^Q_1 = \cdots = g^Q_n=0
$$
where $g_j = x_j\frac{\partial g}{\partial x_j}$ have no common solution on the torus $T$. 

By \cite[Prop.~4.16]{Batyrev (1993)} this is equivalent with the condition that $\bar{Z}_g$ has smooth intersection with all strata $T(Q)$ of the toric variety $\mathbb{P}(P)$, where $Q$ runs over the faces of $P$. 
\begin{lemma}\label{transverse}
Let $U \subset X$ be a toroidal embedding and let $D \subset X$ be a codimension one subvariety which intersects all strata of $X \setminus U$ transversely. Then 
$D$ has toroidal singularities and $U\cap D  \hookrightarrow D$ is a toroidal embedding.
\end{lemma}
{\em Proof }
See \cite[Sect.~13.2]{Danilov (1978)}.

\subsection{Newton modification}
Consider the inclusion $\Delta \subset \sigma^\vee$ for a convenient function $f$. Though $\Delta$ is not a polytope, it corresponds to a toric modification of $X_\sigma$ as follows (see \cite[Sect.~2.1]{Danilov (1979)}. For each face $\Gamma$ of $\Delta$ we let $\sigma_\Gamma \subset \sigma$ be the cone spanned by all $u'-u$ with $u \in \Gamma$ and $u'\in \Delta$. Then the $\sigma_\Gamma$ form a subdivision $\Sigma_\Delta$  of  $\sigma$, called the \emph{polar fan} of $\Delta$. 
Hence the corresponding map $\pi: \mathbb{P}_{\Delta} := X_{\Sigma_\Delta} \to X_\sigma$
is a toric modification. It has been considered by Varchenko \cite{Varchenko (1976)} and Danilov \cite{Danilov (1979)}. Note that $\pi:  \mathbb{P}_{\Delta} \setminus \pi^{-1}(x) \to X_\sigma \setminus \{x\}$ is an isomorphism. 
\begin{theorem}\label{temb}
Let $f \in A_{\sigma}$ be convenient and Newton nondegenerate. Then there exists a closed subset $B \in X_\sigma$ with $x \not\in B$ such that  the inclusion 
\begin{equation}
T\setminus (f^{-1}(0) \cup B) \hookrightarrow  \mathbb{P}_{\Delta}\setminus B
\end{equation}
is a toroidal embedding.
\end{theorem}
{\em Proof }
Let $\theta_1,\ldots,\theta_n$ be a basis of invariant vector fields on $T$ and let $Z_i \subset T$ be given by $\theta_i(f) \neq 0$. Further, let $\bar{Z_i}$ be the closure of $Z_i$ in $ \mathbb{P}_{\Delta}$. Then the fact that $f$ is Newton nondegenerate implies that 
\begin{equation}
\bigcap_{i=1}^n \bar{Z_i} \cap \pi^{-1}(x) = \emptyset
\end{equation}
by \cite[Prop. 4.3]{Batyrev (1993)}. We take $B =\bigcap_{i=1}^n \bar{Z_i}$. Then $B$ is closed in  $ \mathbb{P}_{\Delta}$.
As $T \hookrightarrow  \mathbb{P}_{\Delta}$ is a torus embedding, evidently $T \setminus B \hookrightarrow  \mathbb{P}_{\Delta} \setminus B$ is a toroidal embedding. 

Consider the toroidal embedding $(T\setminus B) \times \mathbb{G}_m \hookrightarrow  ( \mathbb{P}_{\Delta} \setminus B) \times \mathbb{A}^1_k$. The graph $D$ of $f\pi: \mathbb{P}_{\Delta} \setminus B \to  \mathbb{A}^1_k$ satisfies the requirements of Lemma \ref{transverse}. Hence $D \cap( (T\setminus B) \times \mathbb{G}_m )\hookrightarrow D \cap ( ( \mathbb{P}_{\Delta }\setminus B) \times \mathbb{A}^1_k)$ is a toroidal embedding. Pulling this back via the embedding of the graph we obtain the result. 

\subsection{Motivic Milnor fibre}
It is our aim to give an explicit formula for the motivic Milnor fibre $\psi_{f,x}$ for $f \in A_\sigma$ Newton nondegenerate. We use the Newton modification $\pi:  \mathbb{P}_{\Delta}\to X_\sigma$ and observe that 
Theorems (\ref{temb}) and (\ref{psiu}) imply that 
\begin{equation}\label{motnearby}
\psi_f =  \sum_s [D_s^0](1-\mathbb{L})^{\sharp J -1} \in K_0(\mathrm{Var}_k^{\hat{\mu}} )
\end{equation}
where for each stratum $s$ of $(f\pi)^{-1}(0)$ we let $D_s^0$ denote the corresponding subvariety of $\tilde{X}_{\Sigma^+}$ with its $\hat{\mu}$-action. We conclude that $\psi_{f,x}$ is given by the same formula, but where $s$ runs only over the strata contained in $\pi^{-1}(x)$. The only remaining difficulty is to describe $D_s^0$ in simple terms. This however has essentially been done already in  \cite[Sect.~3]{Danilov (1979)}. We briefly recall the result. 

First note that each stratum of $\pi^{-1}(x)$ in $ \mathbb{P}_{\Delta}$ is isomorphic to an algebraic torus $T_\Gamma$. In fact, the closure of such a stratum is of the form $\mathbb{P}_{\Gamma}$ for a compact face $\Gamma$ of $\Delta$. The strata of  $\pi^{-1}(x)$ in the toroidal embedding (\ref{temb}) are then of two types: either the zero set $Z_s$ of $f^{\Gamma}$ in $T_\Gamma$ or $T_\Gamma \setminus Z_{f^\Gamma}$. It remains to describe the strata of the variety $\tilde{X}_{\Delta}$.
\begin{lemma}
Let $\tilde{\Delta}$ be the Newton polyhedron of the function $F = f(z) - t^e \in A_\sigma[t]$. Then $F$ is Newton non-degenerate.
 Moreover, the variety $\tilde{X}_{\Delta}$ is isomorphic to the strict transform of the divisor of the function $F$ in $\mathbb{P}_{\tilde{\Delta}}$. The $\mu_e$-action on these varieties is induced by multiplication of the coordinate $t$ by roots of unity. 
\end{lemma}

\begin{lemma}
For any compact face $\Gamma$ of $\Delta$ let $\hat{\Gamma}$ be the convex hull of $\Gamma$ and $\{0\}$. Then the closures of strata in $\mathbb{P}_{\tilde{\Delta}}$ contained in $(\pi\rho)^{-1}(x)$ are exactly the $\mathbb{P}_{\hat{\Gamma}}$ and the corresponding strata of  $\tilde{X}_{\Delta}$ are either the $U'_\Gamma:=Z_{f^\Gamma } \subset T_\Gamma$ or the $U_\Gamma:=Z_{f^\Gamma-1 } \subset T_{\hat{\Gamma}}$.
\end{lemma}
 Note that for each face $\Gamma$ of $\Delta$ there exists a unique face $\tau$ of $\sigma^\vee$ with $\Gamma^0 \subset \tau^0$ (recall that $\tau^0$ is the relative interior of $\tau$). Let $I_\tau$ denote the set of these compact faces of $\Delta$. The toric stratification of $X_\sigma \setminus B$ induces a decomposition $\psi_f = \sum \psi_{f^\tau}$ and $U_\Gamma$ and $U'_\Gamma$ contribute to $\psi_{f^\tau,x}$ exactly when $\Gamma \in I_\tau$. Let $c_\Gamma = \dim \tau - \dim \Gamma$ for $\Gamma \in I_\tau$. Then we find by Theorem \ref{psiu}:
\begin{theorem}\label{motmil}
\begin{equation}
\psi_{f,x} = \sum_\tau \sum_{\Gamma \in I_\tau} \left((1-\mathbb{L})^{c_\Gamma-1}[U_\Gamma] + (1-\mathbb{L})^{c_\Gamma}[U'_\Gamma]\right).
\end{equation}
\end{theorem}

\section{Example: weighted homogeneous Laurent polynomials}
\begin{definition} Let $T = T_N$ be an $n$-dimensional algebraic torus over $k$ with coordinate ring $k[M]$. 
 A Laurent polynomial $f \in k[M]$  is called \emph{weighted homogeneous} if there exists a  nonzero linear form $\ell \in \mathrm{Hom}(M,\mathbb{Q})$ such that the support of $f$ is a subset of $\ell^{-1}(1)$. 
\end{definition}
Let $f \in k[M]$ be weighted homogeneous with respect to the linear form $\ell$. 
 This linear form determines a positive integer $e$ by setting $\ell(M) = \frac 1e \mathbb{Z}$. It also determines 
 an action  $\gamma^\ast$ of   $\mu_e$ on $k[M]$ and hence a dual action $\gamma$  on  $T_N$ by 
$$\gamma(\zeta)^\ast(X^m) = \zeta^{e\ell(m)} X^m.$$
Here we write $X^m$ for the monomial corresponding to $m \in M$ in $k[M]$. 
Note that $f$ is invariant under this action. 

Let $\Gamma_f$ denote the convex hull $\mathrm{Conv}(\mathrm{Supp}(f))$ in $M\otimes \mathbb{R}$ of the support of $f$. Then $f$ is weighted homogeneous if and only if $\Gamma_f$ is contained in an affine hyperplane not passing through $0$. We will assume that $\dim \Gamma_f = n-1$ and that $f$ is $\Gamma_f$-regular. In that case the linear form $\ell$ is uniquely determined by $f$. 

We have an $n$-dimensional cone $\sigma^\vee \subset M_\mathbb{R}$ given by 
$$\sigma^\vee = \bigcup_{k\geq 0} \mathrm{Conv}(0,k\Delta_f).$$
We let $A_\sigma$ denote the semigroup algebra $k[M \cap \sigma^\vee]$ and $X_\sigma := \mathrm{Spec}(A_\sigma)$.  
Let $\hat{\Delta}_f := \mathrm{Conv}(O,\Delta_f)$. It is a compact $n$-dimensional  polyhedron and corresponds to a projective $n$-dimensional toric variety $\mathbb{P}_{\hat{\Delta}_f}$, which contains $X_\sigma$ as a dense Zariski-open subset. Clearly, the function $f-1$ is $\hat{\Delta}$-regular and invariant under the action of $\mu_e$. We put $V$ the zero set of $f-1$ in $\mathbb{P}_{\hat{\Delta}}$ with its $\mu_e$-action and  $V_\infty$ the zero set of $V$ in $\mathbb{P}_\Delta$. Then 
\begin{equation}
\psi_f = [V] - [V_\infty]\ \mbox{ and } \psi_{f,x} = [V] - \mathbb{L}[V_\infty] .
\end{equation}
Indeed, the first formula follows from Theorem \ref{motnearby} as follows. To each face of $\sigma^\vee$ correspond only three strata of $(f\pi)^{-1}(0)$, namely the intersection of the exceptional divisor with the strict transform of $f^{-1}(0)$ and its complement in these two divisors. Adding corresponding terms for all faces of $\sigma^\vee$ and passing to the $e$-fold ramified cover we obtain $V_\infty$, $f^{-1}(0)\setminus \{0\}$ and  $V\setminus V_\infty$ respectively. Hence 
\begin{equation}
\psi_f = [V\setminus V_\infty] + [f^{-1}(0)\setminus \{0\}] + (1-\mathbb{L})[V_\infty] =  [V]- [V_\infty]
\end{equation}
because $ f^{-1}(0)\setminus \{0\}$ is a $\mathbb{G}_m$- bundle over $V_\infty$, hence $$ [f^{-1}(0)\setminus \{0\}]  = (\mathbb{L}-1)[V_\infty].$$ 
In the formula for $\psi_{f,x}$ we have to subtract the term  $[ f^{-1}(0)\setminus \{0\}]$ from this, so 
$$\psi_{f,x} =  [V]+  (1-\mathbb{L})[V_\infty] = [V] - \mathbb{L}[V_\infty].$$

\end{document}